\documentclass[preprint,12pt]{elsarticle}

    \makeatletter
    \def\ps@pprintTitle{%
      \let\@oddhead\@empty
      \let\@evenhead\@empty
      \def\@oddfoot{\reset@font\hfil\thepage\hfil}
      \let\@evenfoot\@oddfoot
    }
    \makeatother

\usepackage{url}

\usepackage{graphicx}
\usepackage{epstopdf} 

\usepackage{amssymb}
\usepackage{amsmath}
\usepackage{textcomp} 

\begin{document}

\begin{frontmatter}


\title{Enumeration of nonisomorphic Hamiltonian cycles on square grid graphs}

\author{Ed Wynn}

\address{175 Edmund Road, Sheffield S2 4EG, U.K.}
\ead{ed.wynn@zoho.com}

\begin{abstract}
The enumeration of Hamiltonian cycles on $2n\times 2n$ grids of nodes is a longstanding problem in combinatorics.  Previous work has concentrated on counting all cycles.  The current work enumerates nonisomorphic cycles -- that is, the number of isomorphism classes (up to all symmetry operations of the square).  It is shown that the matrix method used previously can be modified to count cycles with all combinations of reflective and 180\textdegree{ }rotational symmetry.  Cycles with 90\textdegree{ }rotational symmetry were counted by a direct search, using a modification of Knuth's Dancing Links algorithm.  From these counts, the numbers of nonisomorphic cycles were calculated for $n\leq10$.
\end{abstract}


\end{frontmatter}


\section{Introduction}
\label{introduction}
We are interested in enumeration of Hamiltonian cycles on $P_{2n} \times P_{2n}$ -- that is, on square grid graphs of $2n \times 2n$ nodes.  The varied interest in this problem can be illustrated by the many names attached to it: as well as Hamiltonian cycles or circuits, there are rook's tours and compact self-avoiding walks.  Since Thompson \cite{Thompson1977} stated the problem and provided elementary results, several authors have tackled this problem.  For example, Myers \cite{Myers1981} converted the problem into enumeration of the skeleton graphs inside the cycle; this was considered to be a simpler problem.  In 1994, Harris Kwong and Rogers \cite{HarrisKwongRogers1994} reported a recent ``flurry of interest'' in the problem. These and other researchers \cite{Kreweras1992,KloczkowskiJernigan1998} use a matrix method, one version of which is discussed in Section~\ref{matrixmethod}.  This method was earlier applied to self-avoiding walks (SAWs) on an infinite square lattice \cite{Enting1980}.  It is still used in that area \cite{ConwayEntingGuttmann1993}.  Jensen \cite{Jensen2004} reports that enumeration of SAWs is `one of the most important and classic combinatorial problems in statistical mechanics'; also it has `traditionally served as a benchmark for both computer performance and algorithm design', and in doing so has had many CPU years applied to it.

\begin{table*}[tb!]
\begin{tabular}{|c| l|}
\hline
Symbol & Symmetries\rule{0pt}{2.5ex}\\
\hline
$u$ & None\rule{0pt}{2.5ex}\\
$v$ & One reflection\\
$w$ & 180\textdegree{ }rotation\\
$x$ & Two reflections and 180\textdegree{ }rotation\\
$y$ & 90\textdegree{ }rotation\\
$z$ & All\\
 \hline
\end{tabular}
\caption{Symbols for numbers of isomorphism classes of Hamiltonian cycles with specified symmetries.  These symmetry descriptions are \emph{exact} -- for example, a cycle with two reflective symmetries is not included in $v$.  The axes of reflection are parallel to the sides of the square.  This list of symmetries is complete.  For example, there is only one cycle with symmetry under diagonal reflection, and this is the trivial cycle in a $2 \times 2$ square, which has all the symmetries of a square.}
\label{symmetrytable}
\end{table*}

These previous efforts have typically counted all reflections and rotations of the cycles, where they are distinct.  For example, the On-Line Encyclopedia of Integer Sequences \cite[{S}equence {A}003763]{OEIS} lists numbers of Hamiltonian cycles on $2n \times 2n$ square grids of nodes up to $n=10$.  The count for $n=11$ is known \cite{Karavaev} and has 71 decimal digits.  By contrast, there are only 4 terms in {S}equence {A}209077, listing the counts of nonisomorphic cycles.  The largest count has 6 decimal digits.  (Here and throughout this paper, a `count of nonisomorphic cycles' is intended to mean the number of equivalence classes of cycles under all symmetry operations of the square.  The title of A209077 calls it a count `reduced for symmetry'.)  There are also variants of Sequence A003763: A143246 counts directed cycles, and therefore contains values twice as large; and A222065 considers $m \times m$ squares.  Hamiltonian cycles are impossible for odd $m$, so A222065 alternates zeros with the values of A143246.  The current work considers only undirected cycles and even $m$.

This paper reports enumeration of nonisomorphic cycles on $2n \times 2n$ square grids of nodes up to $n=10$.  Essentially, the method is to enumerate the symmetric cycles, and to make the appropriate reductions to Sequence A003763.  For most symmetries, the efficient matrix methods can be adapted.  For 90\textdegree{ }rotational symmetry, the enumeration can be done by a direct search on only one quadrant of the square.

\section{Symmetrical cycles}
\label{symmetrics}
Table~\ref{symmetrytable} defines symbols for the numbers of isomorphism classes of Hamiltonian cycles on $2n \times 2n$ grids of nodes. These can be deduced from the counts defined in Table~\ref{counttable}; the matrix equation

\begin{table*}[tb!]
\begin{tabular}{|c| l|}
\hline
Symbol & Symmetries\rule{0pt}{2.5ex}\\
\hline
$A$ & None\rule{0pt}{2.5ex}\\
$B$ & Reflection in a specified axis\\
$C$ & 180\textdegree{ }rotation\\
$D$ & Two reflections and 180\textdegree{ }rotation\\
$E$ & 90\textdegree{ }rotation\\
$F$ & All\\
 \hline
\end{tabular}
\caption{Symbols for counts of Hamiltonian cycles with \emph{at least} the specified symmetries.  The counts are not reduced by symmetry -- for example, $B$ includes two isomorphic copies of each cycle with two reflections.  The axis for $B$ must be horizontal or vertical, not diagonal.}
\label{counttable}
\end{table*}

\begin{equation}\label{matrix1}
\begin{pmatrix} A\\ B\\ C\\ D\\ E\\ F \end{pmatrix} = 
\begin{pmatrix} 8 & 4 & 4 & 2 & 2 & 1\\
                              & 2 &   & 2  &   &  1\\
                              &    & 4 & 2 & 2 & 1\\
                              &    &    & 2  &   & 1\\
                              &    &    &     & 2 & 1\\
                              &    &    &     &    & 1   \end{pmatrix} 
\begin{pmatrix} u\\ v\\ w\\ x\\ y\\ z \end{pmatrix} 
\end{equation}

\noindent can be inverted:

\begin{equation}\label{matrix2}
\begin{pmatrix} u\\ v\\ w\\ x\\ y\\ z \end{pmatrix} = \frac{1}{8}
\begin{pmatrix} 1 &-2 &-1 & 2 &   &  \\
                              & 4 &   &-4  &   &   \\
                              &    & 2 &-2 &-2 & 2\\
                              &    &    & 4 &   &-4\\
                              &    &    &     & 4 & -4\\
                              &    &    &     &    & 8   \end{pmatrix} 
\begin{pmatrix} A\\ B\\ C\\ D\\ E\\ F \end{pmatrix} 
.\end{equation}

Equation~\ref{matrix1} simply expresses how many times each member of an isomorphism class is included in each count.  The coefficients for calculating $B$ depend on the line of reflection having a specified orientation.  So for example, two members from each class in $v$ are included, rotated 180\textdegree{ }relative to each other; two members from each class in $w$ are included, rotated 90\textdegree{ }relative to each other.

Count $A$ is OEIS {S}equence {A}003763.  Section~\ref{matrixmethod} shows that counts $B$, $C$ and $D$ can be generated by small modifications of the matrix method.  Count $F$ is 1 for the $2 \times 2$ square, and 0 for all others.  Count $E$ is considered separately, in Section~\ref{dancinglinks}.  Results are presented in Section~\ref{results}.

\section{The matrix method}
\label{matrixmethod}
The matrix method described here is essentially the same as that of Kloczkowski and Jernigan \cite{KloczkowskiJernigan1998}, with some differences in terminology and implementation.

We define a \emph{connectivity state} (or simply `state') to be a specification of the horizontal edges between two adjacent columns of nodes, together with how these edges are connected in the paths to the left of these columns.  The key observation of the method is that the ways of continuing the paths to the right are independent of how the connectivity state was reached from the left.

\begin{figure}
\centering
\includegraphics[height=5cm]{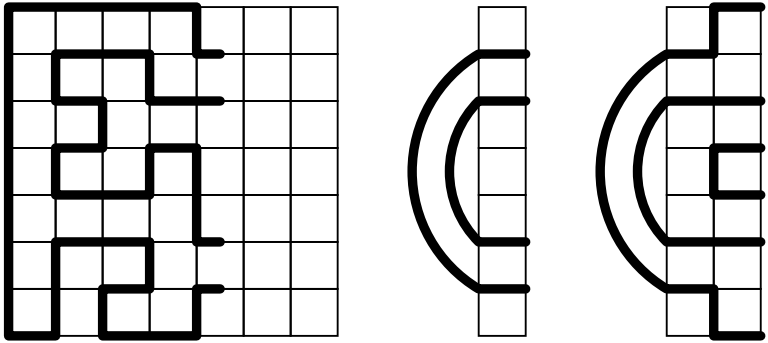}
\caption{The partial cycle in the grid on the left can be represented as the state in the middle; all other ways of reaching the same state will have the same options for continuation.  An example of a continuation is shown on the right: all the nodes in the middle column are visited, and a new state is formed.  In general, it is also possible for a continuation to join links, so long as no loop is formed prematurely.}
\label{fig:state_etc}
\end{figure}

We define a \emph{continuation} $t(s_1;s_2)$ to be a set of vertical edges in the column of nodes to the right of state $s_1$ that ensures that all the column's nodes are visited, and that results in $s_2$ as the next state.  An example is shown in Figure~\ref{fig:state_etc}.  A matrix $\mathbf{T}$ is formed, with each entry $t_{ij}$ equalling 1 or 0 according to whether $t(s_i;s_j)$ exists.  If $\mathbf{\nu}_i$ is a vector of the number of ways of reaching states between columns $i$ and $i+1$, then 
\begin{equation}\label{eq:continuation}
\mathbf{\nu}_{i+1} = \mathbf{T}\mathbf{\nu}_i.
\end{equation}

Some states can be generated by the leftmost column of a square; these are called \emph{starting states}.  (Examples are shown in Tables~\ref{start6counts} and~\ref{start8counts}.)  An equal number of \emph{ending states} can be made into complete Hamiltonian cycles by the rightmost column.  It can be shown that the number of starting states for a $2n\times 2n$ square is the Fibonacci number $F_{2n-1}$.  One state is both starting and ending: this is the state that links the top and bottom nodes.

The count $A$ can therefore be calculated:
\begin{equation}
A = \mathbf{\omega}^T \mathbf{T}^{2n-2} \mathbf{\alpha},
\end{equation}
where $\mathbf{\alpha}$ and $\mathbf{\omega}$ are vectors whose entries are the characteristic functions of starting and ending states respectively.  These vectors are simple to generate; the only remaining problem is how to generate $\mathbf{T}$.  Kloczkowski and Jernigan \cite{KloczkowskiJernigan1998} considered \emph{bond distributions}, which are sets of vertical edges at the right edge of a state.  A different method was used in the current work; see Section~\ref{details}.

An immediate modification of the method is to reject any states or continuations without reflective symmetry in the horizontal axis.  Using overbars to symbolise this modification, the result is count $B$ (defined in Table~\ref{counttable}):
\begin{equation}\label{firstB}
B = \mathbf{\bar{\omega}}^T \mathbf{\bar{T}}^{2n-2} \mathbf{\bar{\alpha}} 
.\end{equation}

Another modification is to apply the unrestricted method only $n-1$ times, to arrive at $\mathbf{\nu}_{n-1}$, the numbers of ways to generate states in the central position.  If and only if a central state consists of only one pair of horizontal edges, the paths that generate that state can be reflected to a create single cycle with reflective symmetry in the vertical axis.  There is only one way to do this in each case, and so
\begin{equation}\label{secondB}
B = \mathbf{\lambda}^T \mathbf{T}^{n-1} \mathbf{\alpha} 
\end{equation}
where $\mathbf{\lambda}$ is a vector containing the characteristic function of single-pair states.  This forms a cross-check for the calculation of $B$ using Equation~\ref{firstB}.

Similarly, some central states will form a single cycle with 180\textdegree{ }rotational symmetry when the paths to the left are copied and rotated to fill the other half of the square.  If $\mathbf{\mu}$ is a vector containing the characteristic function of these states, then
\begin{equation}
C = \mathbf{\mu}^T \mathbf{T}^{n-1} \mathbf{\alpha} 
.\end{equation}
\noindent It is easy to detect single-pair states, to form the vector $\mathbf{\lambda}$ . Various properties can be deduced for the states in vector $\mathbf{\mu}$, but a pragmatic test is simply to rotate and flip the state, and test whether a cycle following these loops does indeed visit all loops and return to the start.

Finally, if the symmetric continuations of Equation~\ref{firstB} are applied to the symmetric  $\mathbf{\bar{\lambda}}$ (or $\mathbf{\bar{\mu}}$, which is identical), then the result has both reflective symmetries and 180\textdegree{ }rotational symmetry:
\begin{equation}
D = \mathbf{\bar{\lambda}}^T \mathbf{\bar{T}}^{n-1} \mathbf{\bar{\alpha}} 
.\end{equation}

Thus, the efficient matrix method can be applied to enumeration of cycles with reflective and 180\textdegree{ }rotational symmetries.  Results are given in Section~\ref{results}. Section~\ref{details} contains details of the implementation of the algorithm and some observations.  The only remaining task is to enumerate 90\textdegree{ }rotational symmetries, and this is considered in Section~\ref{dancinglinks}.

\section{Details and observations from the modified matrix method}
\label{details}

In the current work, continuations were generated from each state by an exhaustive backtracking search.  Starting from each state, each node in the column on the right was considered in turn: a vertical edge and/or a horizontal edge were added to it such the node had degree exactly 2.  It was forbidden for edges to complete a cycle.  When contradictions or complete continuations were reached, the search backtracked -- for example, vertical edges were replaced with horizontal edges.  Thus, the search generated every state that could be reached by a continuation.

Kloczkowski and Jernigan \cite{KloczkowskiJernigan1998} note that some states are impossible because the connections cross, which is not possible in the planar paths.  Also, some states are impossible by parity arguments (similar to those mentioned in Section~\ref{dancinglinks}).  In the current implementation, each iteration considers only those states that have been reached, beginning with the relatively small number of starting states.  Thus, the first iterations need to consider fewer states and are quicker.  Table~\ref{statetable} shows the final numbers of states and continuations that were reached.  It was impractical to store all the continuations, and therefore the approach taken was to store a vector $\mathbf{\nu}_i$, to generate all continuations, and to add each continuation's contribution to $\mathbf{\nu}_{i+1}$ individually.

\begin{table*}[tb]
\begin{tabular}{|c | r r | r r|}
\hline
$n$ & \multicolumn{2}{c|}{Unrestricted} & \multicolumn{2}{c|}{Reflective symmetry\rule{0pt}{2.5ex}}\\
\cline{2-5}
& { }{ }{ }{ }{ }{ }{ }States & { }{ }Continuations & { }{ }{ }{ }{ }{ }{ }States & { }{ }Continuations\rule{0pt}{2.5ex}\\
\hline
1 & 1 & 1 & 1 & 1\rule{0pt}{2.5ex}\\
2 & 6 & 14 & 4 & 6\\
3 & 32 & 162 & 14 & 20\\
4 & 182 & 1966 & 40 & 101\\
5 & 1117 & 25567 & 120 & 327\\
6 & 7280 & 351880 & 320 & 1560\\
7 & 49625 & 5056350 & 946 & 5333\\
8 & 349998 & 75100735 & 2496 & 24727\\
9 & 2535077 & 1144833705 & 7418 & 88422\\
10 & 18758264 & 17821104101 & 19616 & 403552\\
 \hline
\end{tabular}
\caption{Numbers of states and continuations in the matrix method for $2n \times 2n$ grids: in the original, unrestricted method, and when states and continuations are constrained to be symmetric under reflection in the horizontal axis.  The numbers of states in the unrestricted method are in agreement with the values (for $n\le7$) in Table~II of Kloczkowski and Jernigan \cite{KloczkowskiJernigan1998}.}
\label{statetable}
\end{table*}

An efficient form of record-keeping was used to enforce the requirement that continuations must not form loops prematurely.  A record was kept of the \emph{destination} $D(p)$ of each node on the right of the state.  When node $p$ currently has degree 1, for example when nodes are linked by the state at the start of the search, the destination is defined to be the node at the other end of the path.  When $p$ currently has degree 0, then $D(p)$ is defined to be $p$ itself.  When $p$ has degree 2, $D(p)$ is allowed to contain historical information.  In a similar way to the Dancing Links algorithm discussed in Section~\ref{dancinglinks}, this information is used when edges are removed and nodes are released from the middle of a path.  When an edge is added between nodes $p$ and $q$, then the destinations can be updated.  In pseudo-code:
\begin{equation}
\begin{array}{rcl} r & \leftarrow & D(p) \\ s & \leftarrow & D(q) \\ D(r) & \leftarrow & s \\ D(s) & \leftarrow & r  \end{array} 
\label{destinations}
\end{equation}
\noindent In the simplest case, $p$ is linked to another node, $r$, and $q$ is linked to $s$.  When the edge between $p$ and $q$ is added, $r$ and $s$ are now linked and their destinations are updated.   When the edge is later removed, the information in $D(p)$ and $D(q)$ is sufficient to find and restore the old destinations.  The cases where $p$, $q$, $r$ and $s$ are not all distinct need to be detected and treated differently during backtracking.

\begin{figure}
\centering
\includegraphics[width=4.5cm]{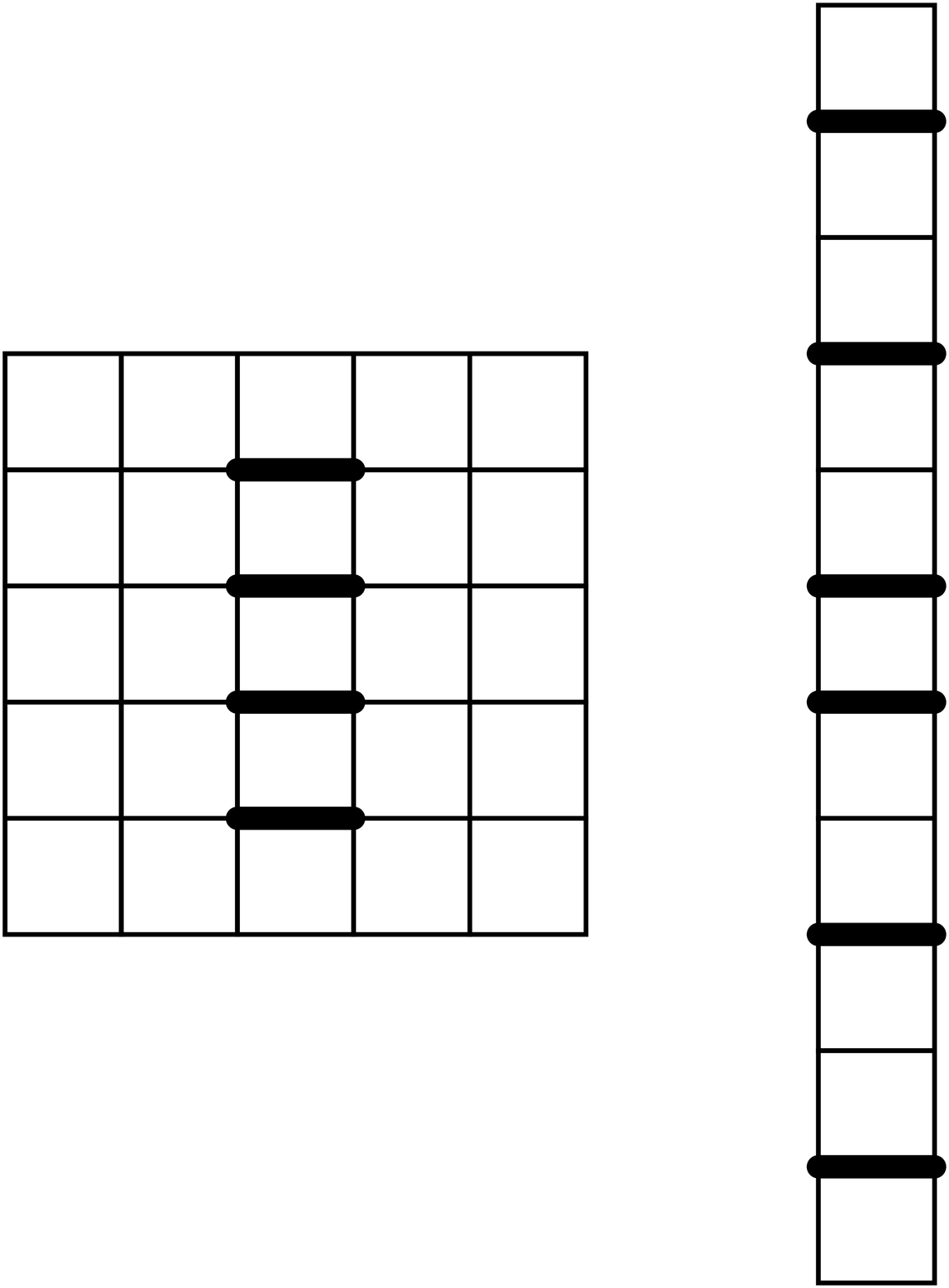}
\caption{Two observations on possible and impossible states.  Left: in the $6\times 6$ grid, there is exactly one way to form a state with the horizontal edges shown in bold (and no others) in the central position.  There is then no way to continue that state to a complete cycle.  Right: in the $12\times 12$ grid, there are many ways to form states with the horizontal edges shown (and no others), but at least 5 continuations are required to form any such state from a starting state.  Therefore, this pattern of horizontal edges can occur \emph{only} in the central position.}
\label{fig:notables}
\end{figure}

The matrix method allows individual states to be considered separately, by starting with special vectors in Equation~\ref{eq:continuation}.  For example, this makes it simple to find the minimum number of continuations to reach a state from starting states, or the minimum number of continuations to reach an ending state.  For example, the grid on the left of Figure~\ref{fig:state_etc} can be completed to a cycle in exactly one way.  Other observations are made in Figure~\ref{fig:notables}.

\begin{figure}
\centering
\includegraphics[width=12cm]{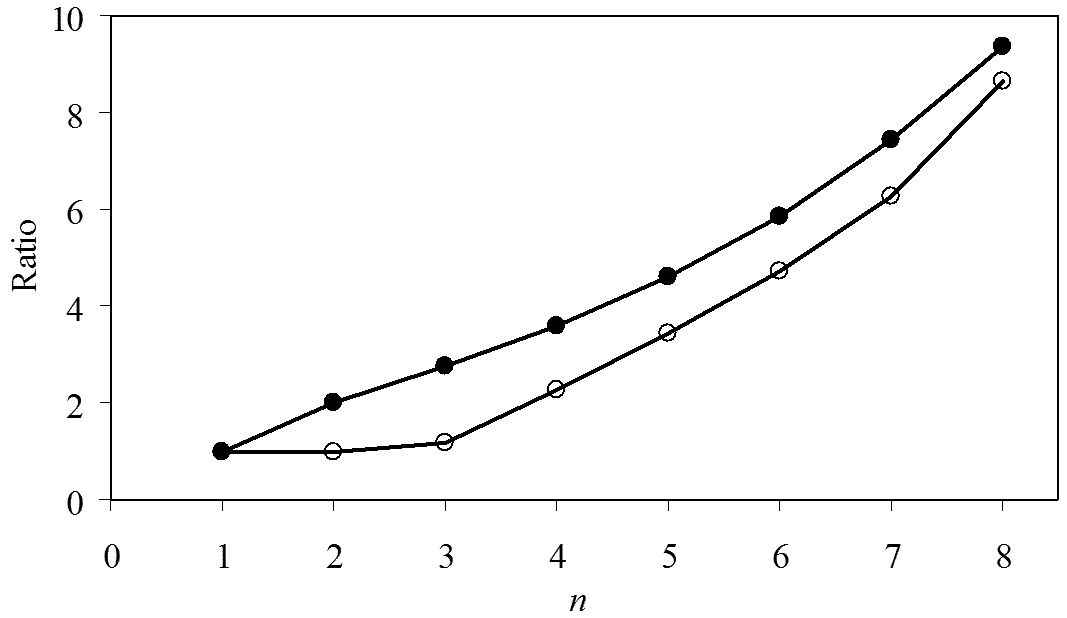}
\caption{Ratios between the counts produced from different starting states on $2n\times 2n$ grids.  Closed symbols: ratio between the single-loop state (the most successul) and the unit-length state.  Open symbols: ratio between the unit-length state and the least successful state.}
\label{fig:ratiograph}
\end{figure}

Tables~\ref{start6counts} and~\ref{start8counts} show the \emph{success rate} of starting states in producing cycles.  The single-loop state is the most successful (for the two sizes shown and all other sizes considered in the current work).  The starting state with unit-length loops (the penultimate state in Table~\ref{start6counts}) is less successful than the single-loop (for $n \geq 2$).  The ratio of counts increase with $n$, as shown in Figure~\ref{fig:ratiograph}.  

Large numbers are reached during the enumerations.  Some previous work (for example \cite{Jensen2004}) used modular arithmetic: the large numbers were stored modulo various integers and the final answers were deduced using the Chinese remainder theorem.  In the current work, the GMP multiple-precision library was used to store and manipulate the values.  This was considered to be more convenient.  From an information-theoretical perspective, there should be no memory penalty in storing the values themselves, and the library allows efficient arithmetic.

\begin{table*}[tb!]
\footnotesize
\begin{tabular}{c c c c c}
\includegraphics[width=1cm]{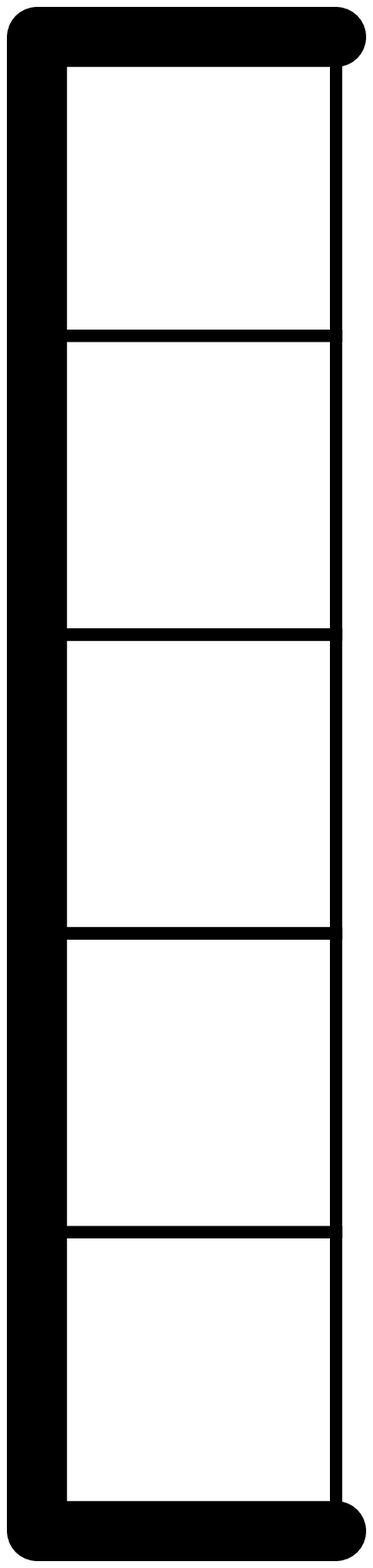} &
\includegraphics[width=1cm]{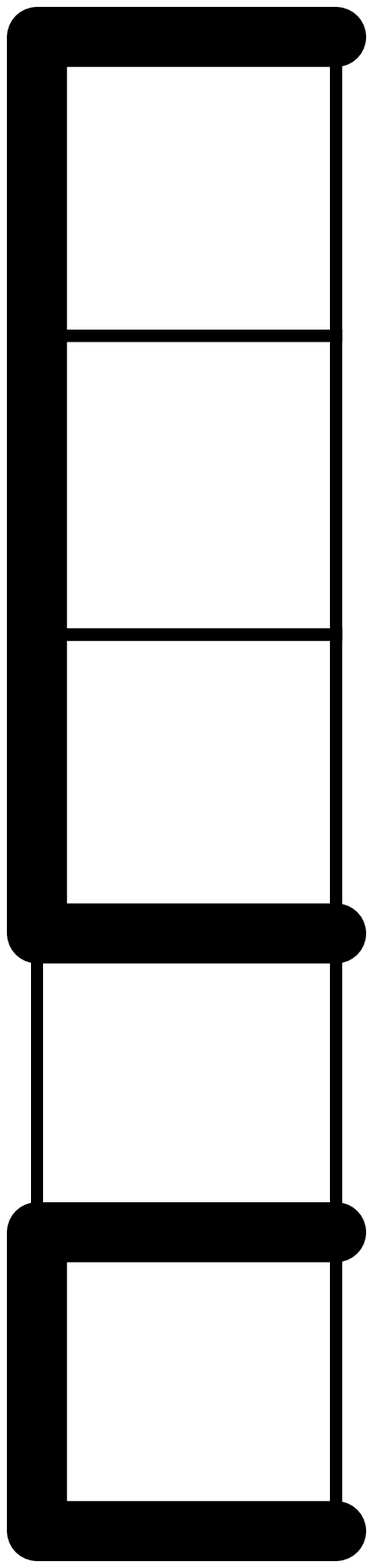} &
\includegraphics[width=1cm]{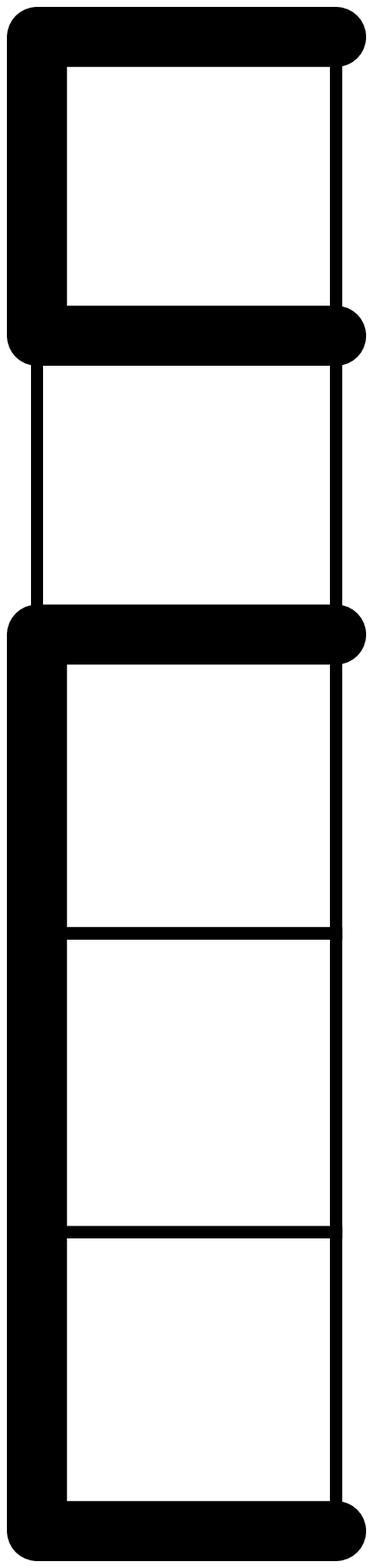} &
\includegraphics[width=1cm]{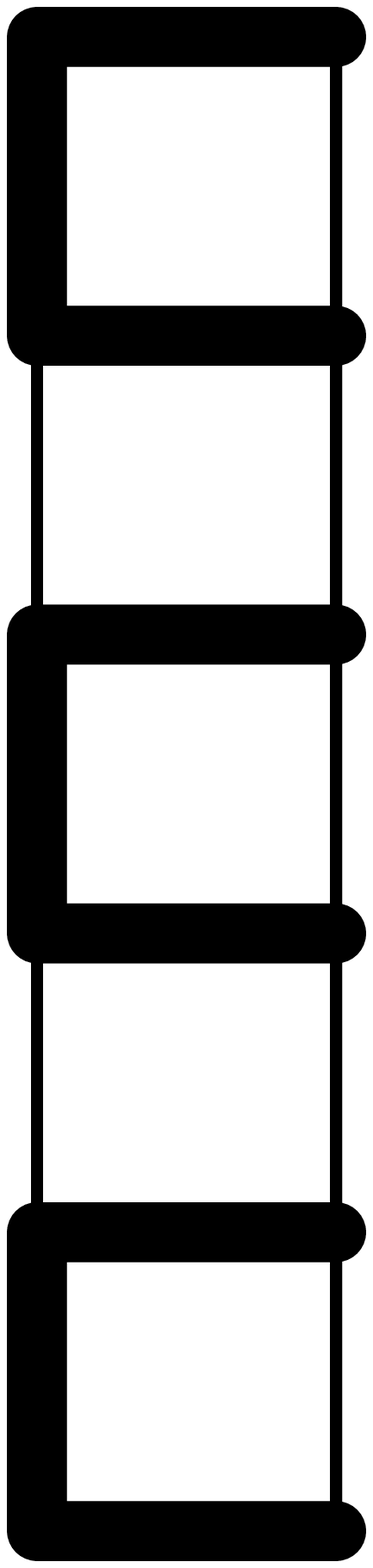} &
\includegraphics[width=1cm]{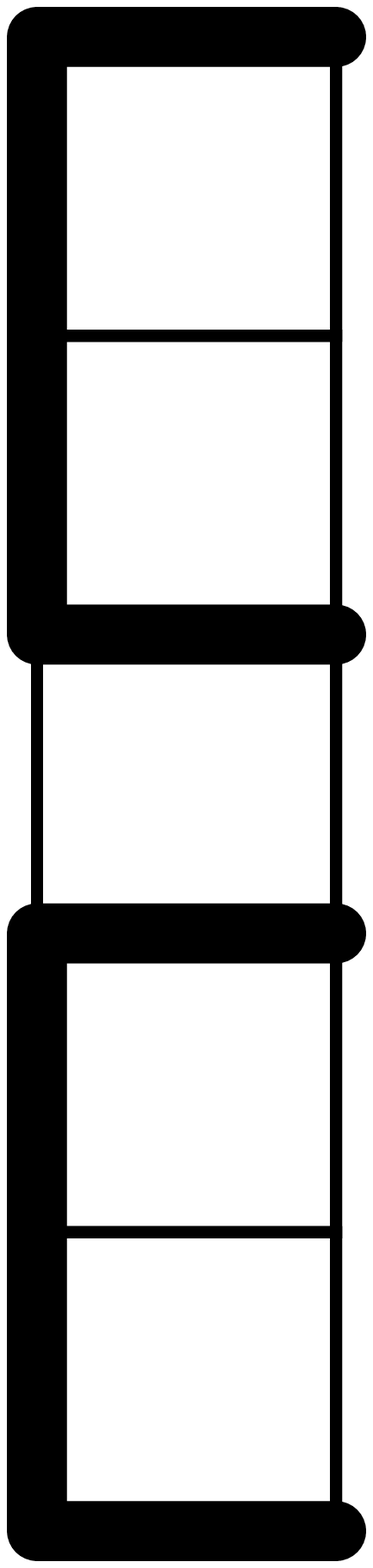} \\
397 & 203 & 203 & 145 & 124\rule{0pt}{2.5ex} \\
\end{tabular}
\caption{Numbers of cycles that can be generated from specified starting states in the $6\times 6$ grid.  The counts are not reduced by symmetry.}
\label{start6counts}
\end{table*}
\begin{table*}[tb!]
\footnotesize
\begin{tabular}{c c c c c c c c c}
\includegraphics[width=1cm]{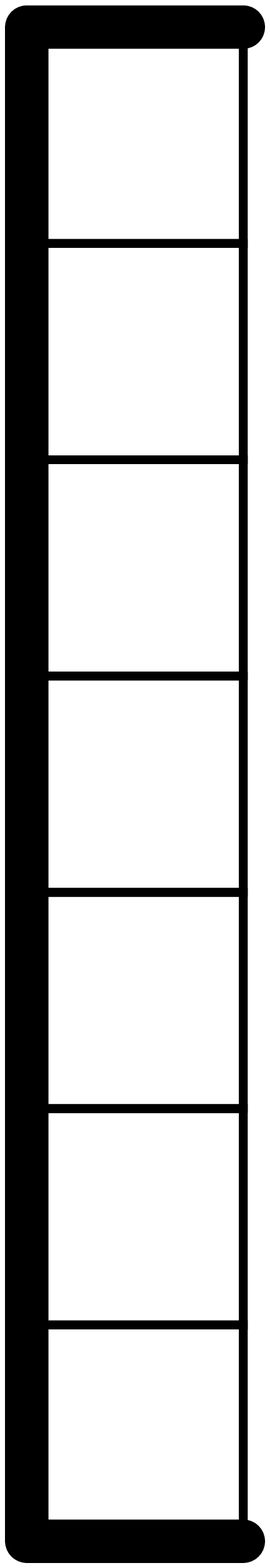} &
\includegraphics[width=1cm]{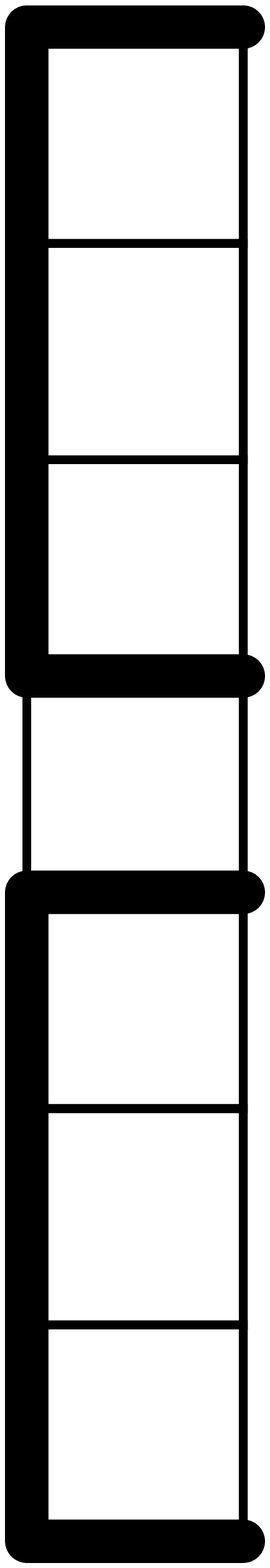} &
\includegraphics[width=1cm]{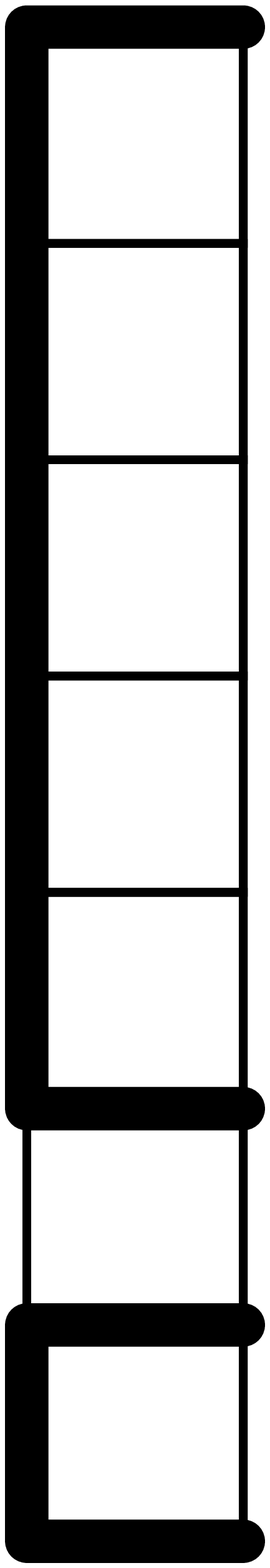} &
\includegraphics[width=1cm]{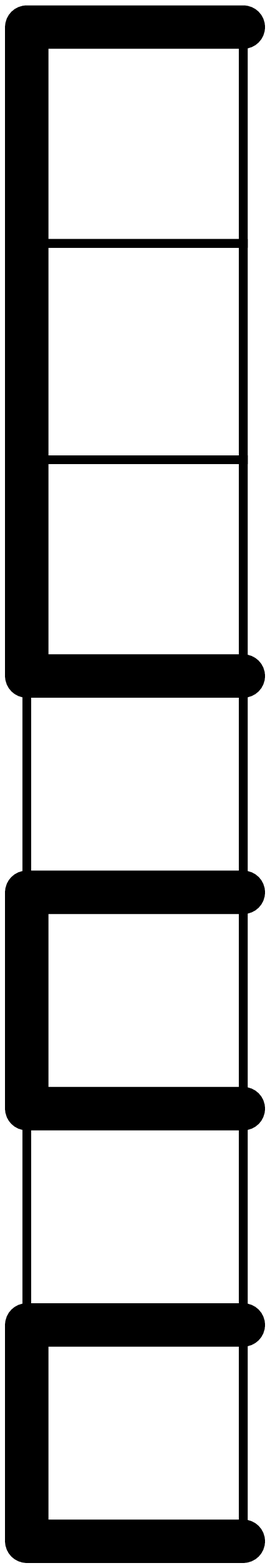} &
\includegraphics[width=1cm]{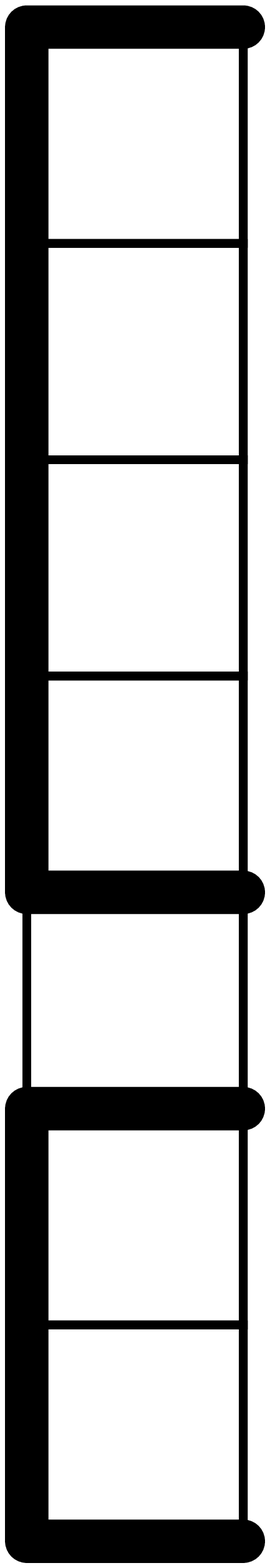} &
\includegraphics[width=1cm]{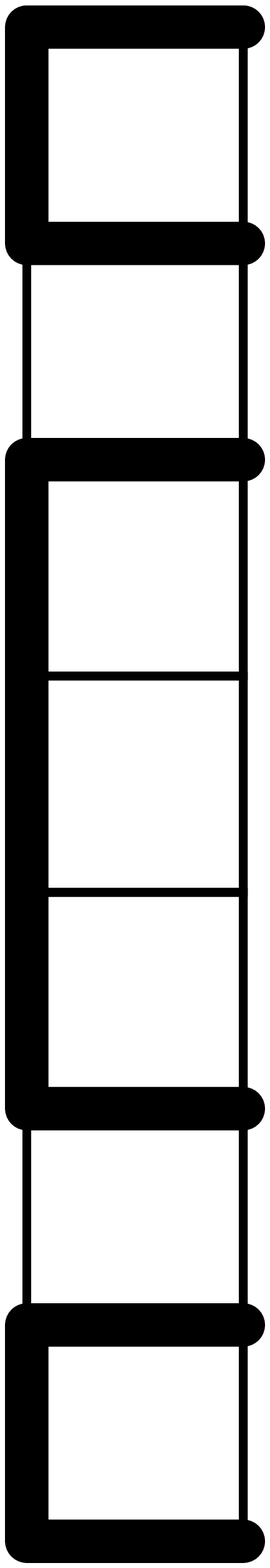} &
\includegraphics[width=1cm]{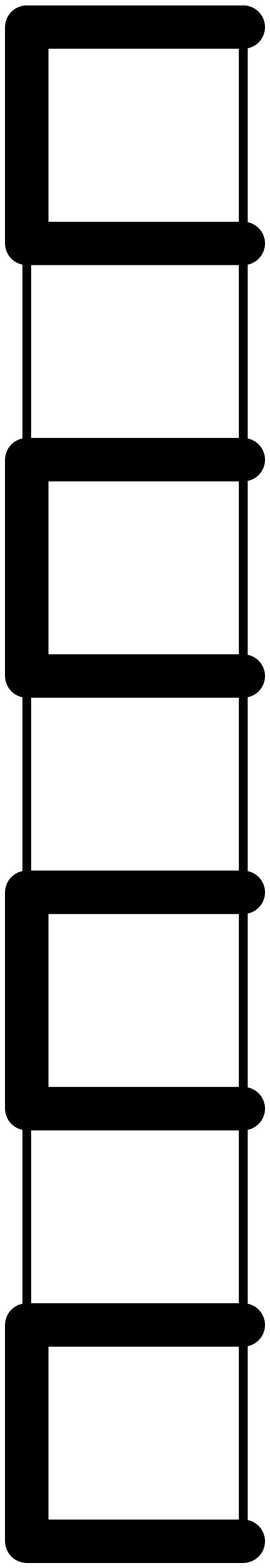} &
\includegraphics[width=1cm]{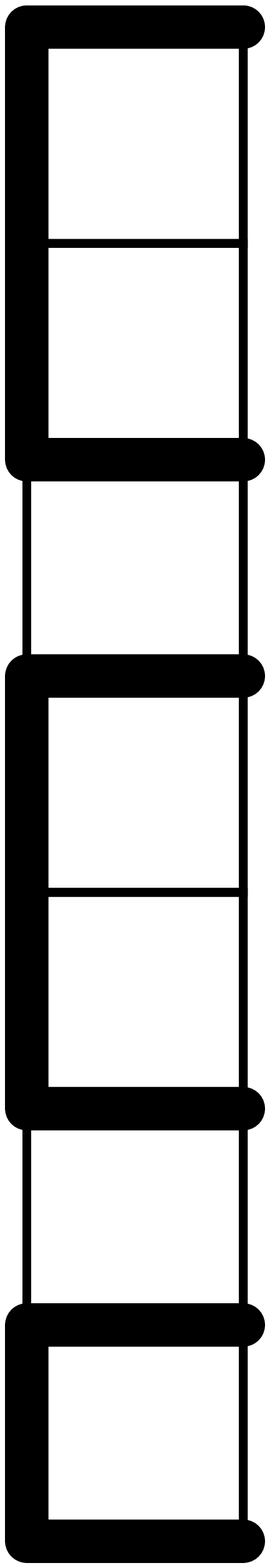} &
\includegraphics[width=1cm]{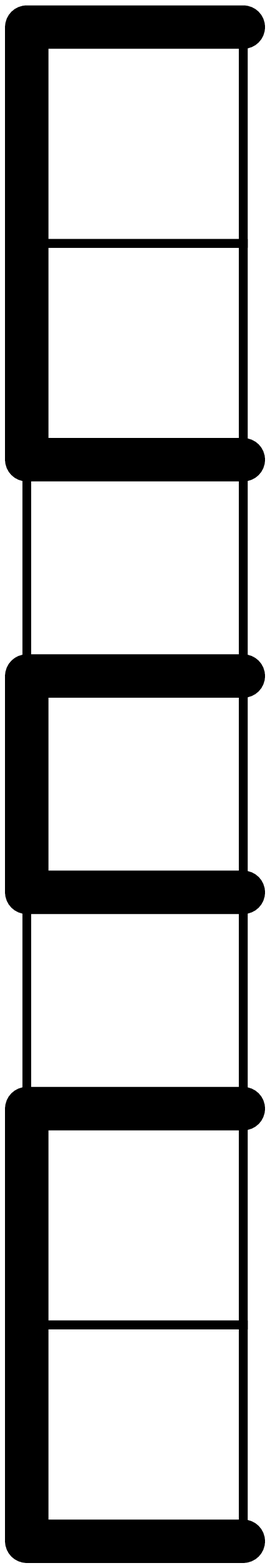} \\
909009 & 510478 & \emph{483465} & \emph{337470} & \emph{322007} & 268967 & 253695 & \emph{149394} & 111755\rule{0pt}{2.5ex} \\
\end{tabular}
\caption{Numbers of cycles that can be generated from specified starting states in the $8\times 8$ grid.  The counts are not reduced by symmetry.  Counts shown in \emph{italics} also apply to the reflected copy of the state.  When these counts are duplicated, the overall total is 4638576, equalling $A_4$ in Table~\ref{counttable}.}
\label{start8counts}
\end{table*}

\section{Cycles with 90\textdegree{ }rotational symmetry}
\label{dancinglinks}
\begin{figure}
\centering
\includegraphics[width=7cm]{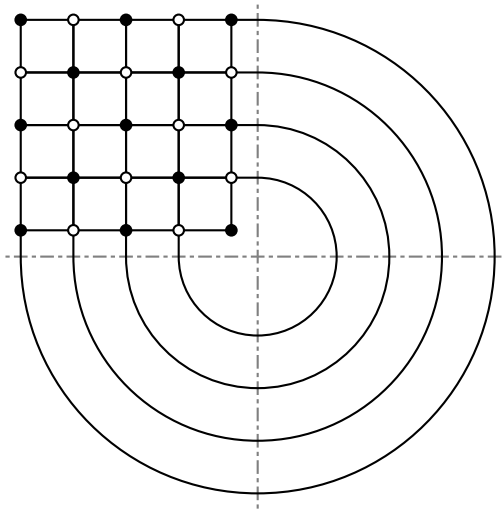}
\caption{A modified $5\times 5$ grid, suitable for searching for $10\times 10$ solutions with 90\textdegree~rotational symmetry.}
\label{fig:grid90}
\end{figure}

A direct search on a modified grid was used to find all cycles with 90\textdegree{ }rotational symmetry.  The modified grid is shown in Figure~\ref{fig:grid90}: for cycles on $2n\times 2n$ grids, a grid of size $n\times n$ is used, with additional edges between the right edge and the bottom edge as shown.  Hamiltonian cycles on the modified grid can be copied to the larger grid, with the additional edges converted to edges between the copies.  An odd number of the additional edges must be used in the cycle.  It is found that this always occurs for odd $n$, and there are no cycles for even $n$.  An outline of a proof is given below.

There is a well-known parity argument proving that Hamiltonian cycles are not possible in $n\times n$ grid for odd $n$.  (For example, see  \cite{KloczkowskiJernigan1998}, where it is also shown by parity arguments that some states are impossible in Hamiltonian cycles.)  If the nodes are coloured alternately as in Figure~\ref{fig:grid90}, then every edge must join nodes of opposite colours.  Therefore, equal numbers of nodes of both colours are present in a cycle.  However, for odd $n$, there is one more node of the `corner' colour, so no cycle can be Hamiltonian.  In the modified grid, by contrast, additional edges link nodes of equal colour.  Therefore, a Hamiltonian cycle is possible if it includes one of these additional edges, linking nodes of `corner' colour.  Further pairs of additional edges, of opposite colours, are also possible.  The odd additional edge permits cycles with odd $n$ but invalidates even $n$.  Therefore, square grids of nodes can have Hamiltonian cycles with 90\textdegree{ }rotational symmetry if and only if they have length $4k+2$ for some $k$.

Knuth's Dancing Links~X (DLX) algorithm \cite{Knuth1999} was modified to perform exhaustive backtracking searches for cycles on the modified grids.  In the terminology of the DLX algorithm, each node was represented by a column, and each edge by a row containing references to its two nodes.  DLX was originally applied to exact cover problems, where rows are selected such that every column is present in exactly one row.  When a row was selected, all columns in that row could be removed from consideration, along with all other rows including those columns.  The `dancing links' make it easier to reinstate these columns when the row is deselected during backtracking.  The modification was to insist that each column must be visited exactly twice.  When a row is selected and a column is visited for the first time, then only the row is removed from consideration.  When the column is visited twice, then that column and all rows that include it are removed.

There is one important detail in this modification: when a column is used as the basis for selecting rows in the search, then a pair of rows should be chosen only once.  Therefore, the search will select a second row only if it follows the first selected row.  Also, when a column is used once as the basis for selecting rows, then it is immediately used again if possible.  These precautions are required to avoid double counting.

This modified method is suitable for direct searches for Hamiltonian cycles on other graphs.  Since every cycle is reached individually, the very large enumerations of Section~\ref{matrixmethod} are not feasible.  However, for the specific application to 90\textdegree{ }rotational symmetry, the reduced search size (on grids of size $n$ rather than $2n$) makes the search manageable.  Results are given in Section~\ref{results}.

\section{Results}
\label{results}

The methods described in earlier sections were implemented on a four-core Intel Core i7-3770 machine with the Windows~8 operating system.  Programs were written in the C language, using the GCC compiler, version 4.8.0, and the Gnu Multiple Precision library, version 5.1.3.  CPU times are quoted from running one process at a time.

For the search using the matrix method, Table~\ref{statetable} shows that the number of continuations increases by a factor of approximately 15 for each increase in $n$.  The complexity of finding each continuation also increases.  Also, the search for all continuations must be done $2n-2$ times, since (in the current work) continuations could not be stored in memory.  Therefore, there were rapid increases in the computer resources required: for grids of $16\times16$,  $18\times18$ and $20\times20$ nodes, the CPU times were 6~seconds, 2~minutes and 16~hours respectively.  The largest of these used approximately 3~GB of RAM, increasing over the course of the run because of the memory requirement for two counts (in vectors $\mathbf{\nu}_i$ and $\mathbf{\nu}_{i+1}$) for every state.

Table~\ref{statetable} shows that the number of continuations is much smaller when reflective symmetry is enforced.  For this variant, even a $22\times22$ grid required only 5~minutes of CPU time and 1~GB of RAM.

For the direct search for cycles with 90\textdegree{ }rotational symmetry, as described in Section~\ref{dancinglinks}, the search for $n=7$ took less than 1~second and the search for $n=9$ took 13.75~hours.  The memory requirements were very small: less than 0.3~MB for $n=9$.

Results are shown in Tables~\ref{rawresults}, \ref{deducedresults} and~\ref{oeisresults}.  The current work agrees with previously published results.  Also, the current work was subject to further cross-checks.  As mentioned earlier, Equations~\ref{firstB} and~\ref{secondB} represent alternative routes to $B$, using substantially different calculations (albeit in the same overall method).  Similarly, the matrix method generates counts of cycles on $2n\times m$ grids of nodes as intermediate results, and these can be compared with the $m\times 2n$ grids' results.  The matrix method could be stopped at only $n-1$ continuations, since there are sufficient values of $A$ already known \cite{Karavaev}.  Instead, the method was run for $2n-1$ continuations to allow cross-checks with those values and the extended Sequence~A222200 \cite{Karavaev}.

The direct search had relatively few cross-checks. The previously published results extend only to $n=4$.  However, the method is only a small modification of a search for cycles in an $n \times n$ grid of nodes.  For that problem, the search gives the correct answer up to $n=8$; this test takes 6~seconds of CPU time.  The test for $n=10$ would be expected to take more than 100000 times longer.

The methods used in this work could be applied to similar problems -- most obviously, non-square rectangular grids of points and cubes in higher dimensions.

\begin{table*}[h!]
\footnotesize
\begin{tabular}{| c | l| l| l| l| l| l|}
\hline
$n$ & 
$A${ }{ }{ }{ }{ }{ }{ }{ }{ }{ }{ }{ }{ } & 
$B${ }{ }{ }{ }{ }{ }{ }{ }{ }{ }{ }{ }{ } & 
$C${ }{ }{ }{ }{ }{ }{ }{ }{ }{ }{ }{ }{ } &
$D${ }{ }{ }{ }{ }{ }{ }{ }{ }{ }{ }{ }{ } & 
$E${ }{ }{ }{ }{ }{ }{ }{ }{ }{ }{ }{ }{ } & 
$F$\rule{0pt}{2.5ex}\\
\hline
1 & 1 & 1 & 1 & 1 & 1 & 1\rule{0pt}{2.5ex}\\
2 & 6 & 4 & 2 & 2 & 0 & 0\\
3 & 1072 & 44 & 28 & 6 & 2 & 0\\
4 & 4638576 & 2828 & 1504 & 40 & 0 & 0\\
5 & \multicolumn{2}{l |} {467260456608}&&&&\\
     &  & 564468 & 520176 & 488 & 204 & 0\\
6 & \multicolumn{2}{l |} {1076226888605605706} &&&&\\
     &  & 754425400 & 696179102 &  13782 & 0 & 0\\
7 & \multicolumn{3}{l |} {56126499620491437281263608}&&&\\
     &  & \multicolumn{2}{l |} {3079904455096} &&&\\
     &  &  & \multicolumn{2}{l |} { 5373177281748}&&\\
     &  &  &  &   757626 & 510718 & 0\\
8 & \multicolumn{3}{l |} {65882516522625836326159786165530572}&&&\\
     &  & \multicolumn{2}{l |} {88444819222239178}&&&\\
     &  &  & \multicolumn{2}{l |} {  166903679914150336} &&\\
     &  &  &  & 95835196  & 0 & 0\\
9 & \multicolumn{4}{l |} {1733926377888966183927790794055670829347983946} &&\\
     &  & \multicolumn{2}{l |} {7685637690960745082050}&&&\\
     &  &  & \multicolumn{2}{l |} {  28636599794306124116062} &&\\
     &  &  &  &24236840344  &31008619522  & 0\\
10 & \multicolumn{5}{l |} {1020460427390768793543026965678152831571073052662428097106} &\\
     & & \multicolumn{3}{l |} {4793315937811919497089287562}&&\\
     &  &  & \multicolumn{3}{l |} {  20262965974179958448766775754}&\\
     &  &  &  &  \multicolumn{2}{l |} {13996574799274} &\\
     &  &  &  &  &  0 & 0\\
%
 \hline
\end{tabular}
\caption{Numbers of $2n \times 2n$ cycles with specified minimum symmetries. Symbols are defined in Table~\ref{counttable}.}
\label{rawresults}
\end{table*}

\begin{table*}[tb]
\footnotesize
\begin{tabular}{|c | l| l| l| l| l| l|}
\hline
$n$ & 
$u${ }{ }{ }{ }{ }{ }{ }{ }{ }{ }{ }{ }{ }{ } & 
$v${ }{ }{ }{ }{ }{ }{ }{ }{ }{ }{ }{ }{ }{ } & 
$w${ }{ }{ }{ }{ }{ }{ }{ }{ }{ }{ }{ }{ }{ } &
$x${ }{ }{ }{ }{ }{ }{ }{ }{ }{ }{ }{ }{ }{ } & 
$y${ }{ }{ }{ }{ }{ }{ }{ }{ }{ }{ }{ }{ }{ } & 
$z$\rule{0pt}{2.5ex}\\
\hline
1 & 0 & 0 & 0 & 0 & 0 & 1\rule{0pt}{2.5ex}\\
2 & 0 & 1 & 0 & 1 & 0 & 0\\
3 & 121 & 19 & 5 & 3 & 1 & 0\\
4 & 578937 & 1394 & 366 & 20 & 0 & 0\\
5 & 58407351059 & 281990 & 129871 & 244 & 102 & 0\\
6 & \multicolumn{2}{l |} {134528360800075421}&&&&\\
     &  & 377205809 &174041330  &  6891 & 0 & 0\\
7 & \multicolumn{3}{l |} {7015812452559988037073365} &&&\\
     &  & \multicolumn{2}{l |} {1539951848735} &&&\\
     &  &  & \multicolumn{2}{l |} { 1343294003351} &&\\
     &  &  &  &   378813 & 255359 & 0\\
8 & \multicolumn{3}{l |} {8235314565328229497795808499821534} &&&\\
     &  & \multicolumn{2}{l |} {44222409563201991} &&&\\
     &  &  & \multicolumn{2}{l |} {  41725919954578785} &&\\
     &  &  &  & 47917598  & 0 & 0\\
9 & \multicolumn{4}{l |} {216740797236120772990968348272561831275923059} &&\\
     &  & \multicolumn{3}{l |} {3842818845468254120853} &&\\
     &  &  & \multicolumn{2}{l |} {  7159149948562719664049} &&\\
     &  &  &  & 12118420172 & 15504309761 & 0\\
10 & \multicolumn{5}{l |} {127557553423846099192878370706037904215158660401579043097} &\\
     & & \multicolumn{3}{l |} {2396657968905952750257244144}&&\\
     &  &  & \multicolumn{3}{l |} {  5065741493544986113047994120} &\\
     &  &  &  &  \multicolumn{2}{l |} {6998287399637} &\\
     &  &  &  &  &  0 & 0\\
 \hline
\end{tabular}
\caption{Numbers of isomorphism classes (under all symmetries of the square) of $2n \times 2n$ cycles with specified  symmetries. Symbols are defined in Table~\ref{symmetrytable}.}
\label{deducedresults}
\end{table*}

\begin{table*}[tb]
\footnotesize
\begin{tabular}{|c | l| l| l|}
\hline
$n$ & 
A209077{ }{ }{ }{ }{ }{ }{ }{ }{ }{ }{ }{ }{ }{ }{ }{ }{ }{ } & 
A227257{ }{ }{ }{ }{ }{ }{ }{ }{ }{ }{ }{ }{ }{ }{ }{ }{ }{ } & 
A227005{ }{ }{ }{ }{ }{ }{ }{ }{ }{ }{ }{ }{ }{ }{ }{ }{ }{ }\rule{0pt}{2.5ex}\\
\hline
1 & 1 & 0 & 0\rule{0pt}{2.5ex} \\
2 & 2 & 1 & 1 \\
3 & 149 & 24 & 4 \\
4 & 580717 & 1760 & 20  \\
5 & 58407763266 & 411861 & 346  \\
6 & 134528361351329451 & 551247139 & 6891 \\
7 & \multicolumn{2}{l |} {7015812452562871283559623} &\\
     &  & 2883245852086 & 634172 \\
8 & \multicolumn{2}{l |} {8235314565328229583744138065519908} &\\
     &  & 85948329517780776 & \multicolumn{1}{l |} {47917598}\\
9 & \multicolumn{3}{l |} {216740797236120772990979350241355889872437894}\\
     &  & \multicolumn{2}{l |} {11001968794030973784902}\\
     &  &  & \multicolumn{1}{l |} {27622729933}\\
10 & \multicolumn{3}{l |} {127557553423846099192878370713500303677609606263171680998}\\
     &  & \multicolumn{2}{l |} {7462399462450938863305238264}\\
     &  &  & \multicolumn{1}{l |} {6998287399637}\\
 \hline
\end{tabular}
\caption{Results for existing sequences in the OEIS for Hamiltonian cycles in $2n \times 2n$ grids. A209077 is the number of all isomorphism classes. A227257 and A227005 are the numbers of classes where the orbits under the symmetry group of the square have 4 elements and 2 elements respectively. Another sequence is A227301, the number of classes where the orbits have 8 elements; this is $u$ in Table~\ref{deducedresults}.}
\label{oeisresults}
\end{table*}

\end{document}